# Constrained recursive kernel density/regression estimation by stochastic quasi-gradient methods

Vladimir Norkin,[1,2] [0000-0003-3255-0405]   Vladimir Kirilyuk[1,3] [0000-0001-6646-5714]

[1] V.M. Glushkov Institute of Cybernetics
of the National Academy of Sciences of Ukraine, Kyiv, 03178 Ukraine &
[2] National Technical University of Ukraine "Igor Sikorsky Kyiv Polytechnic Institute"
`vladimir.norkin@gmail.com`
[3] Center for Evaluating the Activity of Scientific Institutions and Scientific Support for
the Development of the Regions of Ukraine of the National Academy of Sciences of Ukraine

**Abstract.** The paper considers nonparametric kernel density/regression estimation from stochastic optimization point of view. The estimation problem is represented through a family of stochastic optimization problems. Recursive constrained estimators are obtained by application of stochastic (quasi)gradient methods to these problems, classical kernel estimates are derived as particular cases. Accuracy and rate of convergence of the obtained estimates are established, asymptotically optimal estimation procedure parameters are found. The case of moving density/regression is particularly studied.

## 1   Introduction

There is an extensive literature on nonparametric density and regression estimation (see, for example, [5, 15, 19, 26, 27, 36, 37, 47, 49, 52, 54]). Recursive kernel density estimators adjust previously obtained estimates with newcoming observations, they were considered in [2, 3, 16, 17, 18, 25, 38, 39, 42, 43, 53], and recursive kernel regression estimates were studied in [33, 34, 40, 41, 44, 46].

Classical kernel estimators deal with a stable unknown density/regression. They are widely used in spatial statistics for smoothing and representation of spatial data and dependencies (see, for example, [20, 21]). But there are many situations where the desired density/regression may change in time, but currently available information is not sufficient to estimate current density/regression exactly. For example, consumers' spatial density, space demand, consumption or production patterns, obviously change in time, but available information is sparse, costly and comes sequentially in time. Another example gives the recognition of moving patterns on a screen on the basis of brightness measurements at some points. Time limitations may force not to process all available information but only update previously obtained pattern through measurements of randomly chosen points on the screen.

Recursive estimation requires minimum computational efforts and makes minimum demands to computer memory and thus allows updating estimates at a dense greed of points with each new observation within tough time limits. So recursive estimation is the most suitable for on-line identification of time varying systems (see, for example, [23, 35, 50]).

Standard approach to take into account time dependencies is to include time variable (or some functions of it) into the list of regressors as an independent variable. This is possible if behaviouristic aspects of the system under consideration can be ignored and time evolution is deterministic and common for all values of other system parameters. We cannot expect such properties from human governed and spatial systems. Human systems exhibit a free hand and stochastic behaviour. For spatial systems time dependence may be different but correlated for different space points. So in the present paper we give up from the idea to reconstruct parametric or nonparametric time model of the system and only intend to track moving densities/regressions describing the system.

Our approach to nonparametric estimation is close to recursive kernel density/regression estimation but is based on nonstationary stochastic approximation/optimization (see [4, 7, 12, 13, 14, 31]). Similar approach was applied in [9] to a econometric models estimation and in [30] to a kernel density estimation, in the latter cases the nonstationarity is reduced to limit extremal problems, i.e. to a sequence of stochastic optimization problem convergent to the limit problem. We formulate nonparametric estimation problem through a family of some time varying constrained stochastic optimization problems (4), (35), and thus take into account possible changes of the density/regression and a priori knowledge (bounds and other constraints) on density/regression values. To track the solutions of these optimization problems we apply stochastic quasi-gradient (SQG) method [6, 7, 8] and get the so-called SQG-estimates.

We study accuracy and rate of convergence of the obtained estimates for any adjustment and kernel parameters (even not necessarily tending to zero) and any number of observations. In a nonstationary case accuracy estimates are based on a priori estimate of the speed of the density/regression changes. Another essential assumption is that we assume boundedness of the time varying quantities during period of observations. We look at the accuracy of SQG-estimates through three criteria: inaccuracy due to a bad initial approximation, asymptotic deterministic error (bias) and asymptotic stochastic error (variance). On the basis of accuracy estimates we find Pareto optimal values of parameters that give maximal accuracy and rate of convergence. The numerical mechanism behind tracking algorithms is based on properties of sequentially estimated convergent toward zero number sequences considered in Appendix. In stationary and quasi-stationary cases we prove strong consistency of the obtained SQG-estimates by means of stochastic Lyapunov function method. To stabilize estimates we apply Cesàro averaging to SQG-estimates, that is a common technique in stochastic optimization (see [24, 29, 48]) but has not been applied in nonparametric estimation.

In essence, we apply certain stochastic optimization method (stochastic quasi-gradient or SQG-method) to nonparametric estimation problems. SQG-method combines Monte Carlo simulation of a system performance with adaptation of system parameters to get better expected performance. But nonstationary estimation problems arise in stochastic optimization itself, for example, as on-line evaluation of values/gradients of the expectation function (see [10]) or probability/quantile functions (see [22]) along with optimization process. Another example gives the problem of tracking the distribution of a stochastic objective function along with optimization process on the basis of scenario generation and computation of the corresponding objective function values. For this problem the density of the objective function values can be recursively reconstructed by nonparametric kernel density estimation technique developed in the paper. A profile of the density and of the distribution function is the basis for calculation different risk measures characterizing underlying decisions.

Thus the paper's particular contribution to the theory of kernel density/regression estimation is (i) the extension of the theory to nonstationary situation, (ii) the development of more general and flexible recursive constrained estimators, and (iii) the elaboration of practical recommendations for the choice of adjustment and kernel parameters depending on the rate of nonstationarity, limit of observations and the demand to accuracy.

The paper proceeds as follows. In Section 2 we outline the density/regression estimation problems and our approach to their solution, namely nonstationary stochastic optimization problems and stochastic quasi-gradient method. In Section 3 we construct SQG and Cesàro estimates of a density function and study their accuracy and convergence under different assumptions on problem functions and procedure parameters. Section 4 extends the results to recursive regression estimates. Section 5 concludes. Appendix summarizes convergence properties of some recursively estimated number sequences.

## 2    Problem setting and approach

The problem of density estimation is to estimate current time varying density $g^t(x)$, $x \in X \subseteq R^r$, $t = 1, 2, ...$, on the basis of observations $\bar{x}_i$ independently sampled from past densities $g^i(x)$, $1 \leq i \leq t$. Physically densities $g^t(\cdot)$ may be representatives of the family $g(\tau, x)$ depending on a continuous time parameter $\bar{\tau}^t$, i.e. $g^t(x) = g(\tau_t, x)$. It is natural to assume that $g(\tau, \cdot)$ continuously depends on $\tau$, so $\delta_t(x) = |g(\tau_{t+1}, x) - g(\tau_t, x)|$ is small if physical time interval $\rho_t$ is sufficiently small, and one can make $\delta_t(x) \to 0$ with prescribed rate choosing the corresponding rate of convergence $\Delta_t \to 0$.

Beside standard requirements $g^t(\cdot) \geq 0$, $\int_X g^t(x)dx = 1$, the unknown densities may possess other properties like monotonicity, concavity and etc., i.e. densities $g^t$ belong to some class $G$. In particular $G$ may be defined by box constraints $G = \{g(\cdot) \mid g(x) \in [\underline{g}(x), \overline{g}(x)], x \in X\}$, that are decomposable in $x$.

Similarly the problem of regression estimation is to estimate time varying regression function $y^t(x)$, $x \in X \subset R^r$, $t = 0,1,...$, from a family of pare input-output observations $(\overline{x}_i, \overline{y}_i), i = 1,2,...,t$, independently sampled from time varying densities $g^i(x,y), x \in X \subseteq R^r, y \in Y \subseteq R^m$. Assume that $y^t(\cdot)$ belongs to a class of functions $Y$, in particular it may be $Y = \{y(\cdot) \mid y(x) \in Y(x), x \in X\}$. The searched regression functions $y^t(x)$ are expressed through densities $g^t(x,y)$ as conditional expectations

$$y^t(x) = \int_{X \times Y} y g^t(y,x) dx dy / g_1^t(x), \quad g_1^t(x) = \int_X g^t(x,y) dy. \tag{1}$$

Here the problem is to track functions $y^t(\cdot)$ through the sequences of observations $\{(\overline{x}_1, \overline{y}_1), ..., (\overline{x}_t, \overline{y}_t)\}$, $t = 0,1,...$.

In what follows we represent the searched density/regression as a solution of some abstract time varying optimization problem:

$$z_t^* = \arg\min_{z \in Z} \Phi^t(z), \quad t = 0,1,...,$$

and the density/regression estimate is an approximate solution $z_t$ of this problem such that $\|z^t - z_t^*\| \to 0$ as $t \to \infty$. The construction of such sequence $z_t$ is called a problem of nonstationary optimization (see [7, 12, 13, 14]). In our case the difficulty is that functions $\Phi^t(z)$ are not explicitly known. So we also introduce close to $\Phi^t(z)$ auxiliary functions $\Phi_\theta^t(z)$ depending on mollifier parameter $\theta$. Functions $\Phi_\theta^t(z)$ appear to be more tractable, since we can construct statistical estimates $\xi_\theta^t(z)$ for their gradients on the basis of observations. Next we construct a sequence of density/regression estimates through the stochastic quasi-gradient procedure:

$$z^{t+1} = \Pi_Z \left( z^t - \rho_t \xi_{\theta_t}^t(z^t) \right), \quad z^0 \in Z, \quad t = 0,1,...,$$

where $\{\rho_t, \theta_t\}$ concordantly tending to zero number sequences, $Z$ is a set of feasible estimates. We prove that under certain conditions the constructed estimates $z^t$ track true values $z_t^*$. It appears that in an unconstrained stationary case under unchanged parameter $\theta_t \equiv \theta$ and $\rho_t = 1/t$ the sequentially obtained estimates $z^t$ coincide with classical kernel density/regression estimates.

## 3    Recursive kernel density estimation

Assume that searched densities $g^t(x)$ belong to the class

$$G = \{g(\cdot) \mid 0 \leq \underline{g}(x) \leq g(x) \leq \overline{g}(x) \leq \overline{g} \quad \forall x \in X\} \tag{2}$$

(more general classes $G$ will be considered in subsection 3.4). The searched density $g^t(x)$ is the minimizer (over $z$) of the function

$$\Phi^t(x,z) = \frac{1}{2}(z - g^t(x))^2, \quad (3)$$

but this function is not explicitly known. Consider another optimization problem:

$$\Phi_\theta^t(x,z) = \frac{1}{2}\int_X \left[z - \frac{1}{\theta^r}K\left(\frac{\bar{x}-x}{\theta}\right)\right]^2 g^t(\bar{x})d\bar{x} \longrightarrow \min_{z \in [\underline{g}(x),\bar{g}(x)]}, \quad \theta > 0, \quad (4)$$

where $x$ is a fixed point, $\theta$ is a positive (mollifier or window) parameter, $K(\cdot) \geq 0$ is some bounded Borel density on $X$,

$$\int_X K(x)dx = 1, \quad \sup_{x \in X \subseteq R^r} K(x) \leq \overline{K}. \quad (5)$$

One can represent

$$\Phi_\theta^t(x,z) = \frac{1}{2}(z - \tilde{g}_\theta^t(x))^2 + \frac{1}{2}\int_X \frac{1}{\theta^{2r}}K^2\left(\frac{\bar{x}-x}{\theta}\right)g^t(\bar{x})d\bar{x} - \frac{1}{2}(\tilde{g}_\theta^t(x))^2, \quad (6)$$

where a minimizer $\tilde{g}_\theta^t(x)$ has the form

$$\tilde{g}_\theta^t(x) = \frac{1}{\theta^r}\int_X K\left(\frac{\bar{x}-x}{\theta}\right)g^t(\bar{x})d\bar{x}, \quad (7)$$

i.e. function $\tilde{g}_\theta^t(x)$ is a mollified density $g^t(x)$.

Under condition

$$\|x\|^r K(x) \to 0 \quad \text{as} \quad \|x\| \to \infty \quad (8)$$

minimizers $\tilde{g}_\theta^t(x)$ of (4) converge to the minimizer $g^t(x)$ of (3),

$$\lim_{\theta \to 0} \tilde{g}_\theta^t(x) = g^t(x) \quad (9)$$

at continuity points of $g^t(x)$ (see [3, 18]), i.e. functions (3) and (6) are closely connected. Further convergence analysis of $\tilde{g}_\theta^t(x)$ to a discontinuous function $g^t(x)$ in terms of epigraphic convergence is available in [11].

If $g^t(x)$ is globally Hölder continuous at $x$, i.e. for a given $x$

$$|g^t(x) - g^t(y)| \leq H(x)\|x - y\|^{\nu(x)}, \text{ for all } y \in X \quad (\nu(x) > 0), \quad (10)$$

then the following estimate holds true

$$|\tilde{g}_\theta^t(x) - g^t(x)| \leq \int_{R^r} K(y)|g^t(x + \theta y) - g^t(x)|dy \leq A(x)H(x)\theta^{\nu(x)}, \quad (11)$$

where

$$A(x) = \int_{R^r} \|y\|^{\nu(x)} K(y)dy < +\infty \quad (12)$$

is a width characteristic of the kernel. If $H(x) \leq H < +\infty$, $A(x) \leq A < +\infty$, and $\nu(x) \geq \nu > 0$ for all $x \in X$ then convergence in (9) is uniform. Remark that a globally Hölder continuous at $x$ function $g^t(x)$ is continuous at $x$ but may be discontinuous at other points. In what follows we as a rule drop argument $x$ in $A(x)$, $H(x)$, $\nu(x)$.

Empirical approximation of $\tilde{g}_\theta^t(x)$ for a given $t$ equals to a classical kernel density estimate

$$z_\theta^{N_t}(x) = \frac{1}{N_t \theta^r} \sum_{i=1}^{N_t} K\left(\frac{\bar{x}_i - x}{\theta}\right), \tag{13}$$

where $\bar{x}_i$, $i = 1,...,N_t$, are independently sampled from the same density $g^t(\cdot)$.

By the law of large numbers for any given $x$ and $\theta$ with probability one $z_\theta^{N_t}(x) \to \tilde{g}_\theta^t(x)$ as $N_t \to \infty$. If density $K(\cdot)$ is continuous, this convergence is uniform in $x$. The accuracy of estimate (13) beside $N_t$ depends also on window parameter $\theta$, optimal choice of $\theta$ depends on $N_t$.

There may be time and cost limitations to apply empirical estimates (13). The observer may have no time and resources to perform/process large number $N_t$ observations $\bar{x}_i$, $i = 1,...,N_t$, within time interval where the searched density $g^t(x)$ can be considered unchanged. For these reasons it makes sense to utilize previously obtained information (samples) from past close densities $g^i(x), i \leq t$. For example, one could consider empirical estimates

$$z_\theta^t(x) = \frac{1}{t\theta^r} \sum_{i=1}^{t} K\left(\frac{\bar{x}_i - x}{\theta}\right), \tag{14}$$

where $\bar{x}_i$, are independently sampled from different $g^i(\cdot)$, $i = 1,...,t$. The question is whether such estimates trace $g^t(\cdot)$ as $t \to \infty$. In the next subsections we consider similar and more general estimators and prove their consistency.

### 3.1 Point-wise SQG-estimates in nonstationary case

Alternatively to (13), we can estimate $g^t(x)$ by solving nonstationary problems (4) with $\theta \to 0$, $t \to \infty$, for example, by SQG-method (see [6, 7, 8]) and obtain a sequence of SQG-estimates:

$$z^{t+1} = \min\{\bar{g}(x), \max\{\underline{g}(x), (z^t - \rho_t \xi^t)\}\}, \quad z^0 \in [\underline{g}(x), \bar{g}(x)], \quad t = 0,1,..., \tag{15}$$

where stochastic gradients

$$\xi_t = \left(z^t - \frac{1}{\theta_t^r} K\left(\frac{\bar{x}_t - x}{\theta_t}\right)\right) \tag{16}$$

are such that conditional expectation $E\{\xi_t | z^t\} = \nabla_z \Phi_{\theta_t}^t(x, z^t)$, nonnegative smoothing parameters $\theta_t$ and adjustment coefficients $\rho_t$ are measurable with respect to $\sigma$-algebra $\sigma\{\bar{x}_0, \bar{x}_1,..., \bar{x}_{t-1}\}$ generated by observations $\sigma\{\bar{x}_0, \bar{x}_1,..., \bar{x}_{t-1}\}$ defined on a joint probability space $(\Omega, \Sigma, P)$, $E$ denotes mathematical expectation over measure $P$. Notice that empirical estimate (13) can be obtained by $N_t$ steps of sequential procedure (15) with $\rho_i = 1/t$, $\theta_t \equiv \theta$ applied to unconstrained problem (4) for a given $t$.

Procedure (15) stems from the theory of nonstationary stochastic optimization (see [7, 12, 13, 14]). For the stationary case $g^t(x) \equiv g(x)$ similar recursive estimation procedures but without projection on interval $[\underline{g}(x), \overline{g}(x)]$ and with deterministic convergent to zero sequences $\{\theta_t, \rho_t\}$ were considered in [2, 3, 18, 53].

**Remark.** If it is possible to make several independent observations $\{\overline{x}_{t1},...,\overline{x}_{tN_t}\}$ from the same density $g^t(\cdot)$ then one can use in (15) another stochastic gradient

$$\xi_t = \left( z^t - \frac{1}{N_t \theta_t^r} \sum_{i=1}^{N_t} K\left( \frac{\overline{x}_{ti} - x}{\theta_t} \right) \right)$$

having a smaller variance than (16).

Simce $\mathbf{E}\{\xi_t \mid z^t\} = \nabla_z \Phi_{\theta_t}^t(x,z) = (z^t - \tilde{g}_{\theta_t}^t(x))$ approximates the gradient $\nabla_z \Phi^t(x,z) = (z^t - g^t(x))$ of function $\Phi^t(x,\cdot)$, then $\xi_t$ is a stochastic quasi-gradient of $\Phi^t(x,\cdot)$ at $z^t$ in terminology of [6, 7, 8]. So the corresponding estimates (15) are called SQG-estimates.

Next theorems establish rate of tracking values $g^t(x)$ by estimators $z^t$ under different assumptions on the speed of density changes.

**Theorem 1** (*accuracy of SQG-estimates for small constant adjustment and window parameters*). Assume that

(i) for all $t$ the tracked densities $g^t(x) \in G$ and their changes satisfy $|g^{t+1}(x) - g^t(x)| \leq \delta(x)$, $\delta(x) > 0$;

(ii) kernel $K(\cdot)$ is a bounded Borel density satisfying (5), (12);

(iii) sequence of estimates $\{z^t\}$ is obtained by (15), where adjustment coefficients $\rho_t$ and window parameters $\theta_t$ are constant, $0 < \rho_t \equiv \rho \leq 1/2$, $\theta_t \equiv \theta$. Then at any point $x$, where all densities $g^t(\cdot)$ are globally Hölder continuous (see (10)), the tracking error satisfies

$$(z^{t+1} - g^{t+1}(x))^2 \leq 4\overline{g}\delta(x)\rho^{-1} + 2\overline{g}AH\theta^v + \overline{g}^2\rho + \overline{K}^2\rho\theta^{-2r}$$
$$+ (1 - 2\rho)^t (z^0 - g^0(x))^2 + \zeta_t,$$

where a stochastic term $\zeta_t$ has zero mean $\mathbf{E}\zeta_t = 0$ and a variance $D\zeta_t \leq 2\rho\overline{g}^3(\overline{g} + \overline{K}\theta^{-r})$.

**Comment 1** (*risk based error minimization*). The error depends on the accuracy of the initial approximation $z^0$, number of steps $t$, the speed of changes $\delta$ and the choice of adjustment parameter $\rho$ and window size $\theta$. The initial error $|z^0 - g^0(x)|^2$ is quickly suppressed with the course of iterations. The error is basically determined by the speed of changes $\delta$, the rest of error terms has to be adapted to this speed by the optimal choice of $\rho$ and $\theta$. The error contains a deterministic term

$$\varphi(\delta, \rho, \theta) = 4\overline{g}\delta\rho^{-1} + 2\overline{g}AH\theta^v + \overline{g}^2\rho + \overline{K}^2\rho\theta^{-2r}$$

and a stochastic one $\zeta_t$, which can be characterized by some risk measure like a variance $D\zeta_t$ or a standard deviation $\sqrt{D\zeta_t}$. So the problem of error minimization over $\rho$, $\theta$ is three criteria, and one can find only Pareto optimal $\rho$, $\theta$. For example, for known constants $\bar{g}$, $\nu$, $H$, $A$, $\bar{K}$ the choice can be made by minimization of an error function $\varphi$ over $\rho \in [0, 1/2]$, $\theta > 0$. If these constants are not known then a reasonable choice of $\rho$ and $\theta$ can be made assuming that $\delta > 0$ and hence $\rho > 0$, $\theta > 0$ are close to zero, i.e. consider the following error function

$$\varphi(\delta, \rho, \theta) = \max\{4\bar{g}, 2\bar{g}AH, \bar{g}^2, \bar{K}^2\} \cdot \left(\delta\rho^{-1} + \theta^\nu + \rho + \rho\theta^{-2r}\right).$$

Taking $\rho = \delta^\alpha$, $\theta = \delta^\beta$ we can find that asymptotically optimal $\alpha^* = \frac{2r+\nu}{2(r+\nu)}$, $\beta^* = \frac{1}{2(r+\nu)}$. Then the error function asymptotically equals

$$\varphi \sim \delta^{\frac{\nu}{2(r+\nu)}}.$$

We have to check also a stochastic error described by the standard deviation

$$\sqrt{D\zeta_t} \leq \sqrt{2}\sqrt{\rho}\bar{g}^2\left(1 + \frac{\bar{K}}{\bar{g}}\theta^{-r}\right)^{1/2} = \psi(\rho, \theta).$$

For the selected $\alpha^*$, $\beta^*$ the stochastic error for small $\delta$ has the following asymptotics

$$\psi \sim \delta^{\frac{1}{2} - \frac{\nu}{4(\nu+r)}}.$$

This shows that the choice $\rho = \delta^{\alpha^*}$, $\theta = \delta^{\beta^*}$ is reasonable: in case of bounded changes $|g^{t+1}(x) - g^t(x)| \leq \delta$ estimates $z^t$ overtake $g^t(x)$, and for possible $\nu \in (0, 1]$ this choice minimizes deterministic error and keeps stochastic error sufficiently small.

**Proof of Theorem 1.** Denote for brevity $g^t = g^t(x)$. Estimate

$$\begin{aligned}(z^{t+1} - g^{t+1})^2 &\leq (z^t - g^{t+1} - \rho_t \xi_t)^2 \\ &\leq (z^t - g^{t+1})^2 - 2\rho_t \xi_t (z^t - g^{t+1}) + \rho_t^2(\xi_t)^2 \\ &\leq (z^t - g^t)^2 + 2(z^t - g^t)(g^t - g^{t+1}) + (g^t - g^{t+1})^2 + \rho_t^2(\xi_t)^2 \\ &\quad - 2\rho_t \nabla_z \Phi^t_{\theta_t}(x, z^t) \cdot (z^t - g^{t+1}) - 2\rho_t \left(\xi_t - \nabla_z \Phi^t_{\theta_t}(x, z^t)\right) \cdot (z^t - g^{t+1}) \\ &\leq (z^t - g^t)^2 + 2(z^t - g^t)(g^t - g^{t+1}) + (g^t - g^{t+1})^2 \\ &\quad - 2\rho_t \nabla_z \Phi^t_{\theta_t}(x, z^t) \cdot (z^t - g^t) - 2\rho_t \nabla_z \Phi^t_{\theta_t}(x, z^t) \cdot (z^t - g^{t+1}) \\ &\quad - 2\rho_t \left(\xi_t - \nabla_z \Phi^t_{\theta_t}(x, z^t)\right) \cdot (z^t - g^{t+1}) + \rho_t^2(\xi_t)^2 \\ &\leq (z^t - g^t)^2 + 2(z^t - g^t)(g^t - g^{t+1}) + (g^t - g^{t+1})^2 + \rho_t^2(\xi_t)^2 \\ &\quad - 2\rho_t \nabla_z \Phi^t(x, z^t) \cdot (z^t - g^t) - 2\rho_t \left(\nabla_z \Phi^t_{\theta_t}(x, z^t) - \nabla_z \Phi^t(x, z^t)\right) \cdot (z^t - g^t) \\ &\quad - 2\rho_t \nabla_z \Phi^t_{\theta_t}(x, z^t) \cdot (g^t - g^{t+1}) - 2\rho_t \left(\xi_t - \nabla_z \Phi^t_{\theta_t}(x, z^t)\right) \cdot (z^t - g^{t+1}).\end{aligned}$$

Taking into account that $\nabla_z \Phi^t(x, z) = (z - g^t)$, $\nabla_z \Phi^t_{\theta_t}(x, z) = (z - \tilde{g}^t_{\theta_t})$ and $|z^t| \leq \bar{g}$, $|g^t| \leq \bar{g}$, $|\tilde{g}^t_{\theta_t}(x)| \leq \bar{g}$, from the above inequality we get

$$\begin{aligned}(z^{t+1} - g^{t+1})^2 &\leq (z^t - g^t)^2 + 2(z^t - g^t)(g^t - g^{t+1}) + (g^t - g^{t+1})^2 + \rho_t^2(\xi_t)^2 \\ &\quad - 2\rho_t (z^t - g^t)^2 - 2\rho_t \left(g^t - \tilde{g}^t_{\theta_t}(x)\right)(z^t - g^t)\end{aligned}$$

$$-2\rho_t\left(z^t - \tilde{g}^t_{\theta_t}(x)\right)(g^t - g^{t+1}) - 2\rho_t\left(\xi_t - \nabla_z\Phi^t_{\theta_t}(x, z^t)\right)\cdot(z^t - g^{t+1})$$
$$\leq (z^t - g^t)^2 - 2\rho_t(z^t - g^t)^2$$
$$+ 4\rho_t\overline{g}\,|g^t - \tilde{g}^t_{\theta_t}(x)| + 4\rho_t\overline{g}\,|g^{t+1} - g^t| + 6\overline{g}\,|g^{t+1} - g^t|$$
$$+ \rho_t^2(\xi_t)^2 + 2\rho_t\left(\xi_t - \nabla_z\Phi^t_{\theta_t}(x, z^t)\right)\cdot(g^{t+1} - z^t). \tag{17}$$

From (11) we have $|\tilde{g}^t_{\theta_t}(x) - g^t| \leq AH\theta_t^\nu$. Since $|g^{t+1}(x) - g^t(x)| \leq \delta_t(x)$,

$$(\xi_t)^2 \leq 2\left((z^t)^2 + \theta_t^{-2r}K^2(\theta_t^{-1}(\overline{x}_t - x))\right) \leq 2\left(\overline{g}^2 + \overline{K}^2\theta_t^{-2r}\right),$$
$$\xi_t - \nabla\Phi^t_{\theta_t}(z^t) = \tilde{g}^t_{\theta_t}(x) - \theta_t^{-r}K(\theta_t^{-1}(\overline{x}_t - x)),$$

then for $\rho_t \leq 1/2$

$$(z^{t+1} - g^{t+1})^2 \leq (z^t - g^t)^2 - 2\rho_t(z^t - g^t)^2$$
$$+ 4\rho_t\overline{g}\,|g^t - \tilde{g}^t_{\theta_t}(x)| + 4\rho_t\overline{g}\,|g^{t+1} - g^t| + 6\overline{g}\,|g^{t+1} - g^t|$$
$$+ 2\rho_t^2\left((z^t)^2 + \theta_t^{-2r}K^2(\theta_t^{-1}(\overline{x}_t - x))\right)$$
$$+ 2\rho_t\left(\tilde{g}^t_{\theta_t}(x) - \theta_t^{-r}K(\theta_t^{-1}(\overline{x}_t - x))\right)\cdot(g^{t+1} - z^t)$$
$$\leq (1 - 2\rho_t)(z^t - g^t)^2$$
$$+ 4\rho_t\overline{g}AH\theta_t^\nu + 8\overline{g}\delta_t(x) + 2\rho_t^2\left(\overline{g}^2 + \overline{K}^2\theta_t^{-2r}\right)$$
$$+ 2\rho_t\left(\tilde{g}^t_{\theta_t}(x) - \theta_t^{-r}K(\theta_t^{-1}(\overline{x}_t - x))\right)(g^{t+1} - z^t). \tag{18}$$

Denote $v_t = (z^t - g^t)^2$, $\eta_t = \left(\tilde{g}^t_{\theta_t}(x) - \theta_t^{-r}K(\theta_t^{-1}(\overline{x}_t - x))\right)\cdot(g^{t+1} - z^t)$ and assume $\rho_t \equiv \rho$, $\theta_t \equiv \theta$, $\delta_t(x) \leq \delta$, then (18) becomes

$$v_{t+1} \leq (1 - 2\rho)v_t + 4\overline{g}AH\rho\theta^\nu + 8\overline{g}\delta + 2\overline{g}^2\rho^2 + 2\overline{K}^2\rho^2\theta^{-2r} + 2\rho\eta_t.$$

By Lemma A1 from Appendix we get

$$v_{t+1} \leq v_0(1 - 2\rho)^t + 2\overline{g}AH\theta^\nu + 4\overline{g}\delta\rho^{-1} + \overline{g}^2\rho + \overline{K}^2\rho\theta^{-2r} + \zeta_t, \tag{19}$$

where a stochastic term $\zeta_t = \sum_{i=0}^{t}\rho(1 - q\rho)^{t-i}\eta_i$ has zero mean $\mathbf{E}\zeta_t = 0$ and a variance

$$D\zeta_t = \sum_{i=0}^{t}\rho^2(1 - 2\rho)^{2(t-i)}E\eta_i^2 \leq \frac{1}{2}\rho\sup_i E\eta_i^2,$$
$$\sup_i E\eta_i^2 \leq \mathbf{E}\,|\tilde{g}^t_\theta(x) - \theta^{-r}K(\theta^{-1}(\overline{x}_i - x))|^2 \cdot |g^{t+1} - z^t|^2$$
$$\leq 4\overline{g}^2(\overline{g}^2 + \theta^{-r}\mathbf{E}\theta^{-r}K^2(\theta^{-1}(\overline{x}_t - x))$$
$$\leq 4\overline{g}^3(\overline{g} + \overline{K}\theta^{-r}). \tag{20}$$

From (19), (20) the required estimate follows.

In the next theorem we assume that speed of changes $|g^{t+1}(x) - g^t(x)|$ tends to zero. Physically this means that either observations $\overline{x}_t$ are made more and more fre-quently in the physical time, or densities $g^t(\cdot)$ are quasi-stationary, for example, approach to a stationary value.

**Theorem 2** (*convergence in mean of SQG-estimates when density changes vanish*). Assume that

(i) for all $t$ the tracked densities $g^t(x)$ belong to class $G$ and their changes satisfy $|g^{t+1}(x) - g^t(x)| \leq \delta_t \to 0$;

(ii) kernel $K(\cdot)$ is a bounded Borel density satisfying (5), (12);

(iii) sequence of estimates $\{z^t\}$ is obtained by (15), where adjustment coefficients $\rho_t$ and window parameters $\theta_t$ are deterministic and satisfy conditions

$$\lim_{t\to\infty}\rho_t = \lim_{t\to\infty}\theta_t = 0, \quad \lim_{t\to\infty}\frac{\delta_t}{\rho_t} = \lim_{t\to\infty}\frac{\rho_t}{\theta_t^r} = 0, \quad \sum_{t=0}^{+\infty}\rho_t = +\infty. \qquad (21)$$

Then at any point $x$, where all densities $g^t(\cdot)$ are globally Hölder continuous (see (10)), holds

$$\lim_{t\to\infty}\mathbf{E}(z^t - g^t(x))^2 = 0.$$

**Comment 2.** Conditions (21) indicate how to choose adjustment and mollifier parameters given the estimates $\delta_t$ of the speed of changes $|g^{t+1}(x) - g^t(x)|$. Taking any sequences $\bar{\rho}_t$, $\bar{\theta}_t$ such that

$$\lim_{t\to\infty}\bar{\rho}_t = \lim_{t\to\infty}\bar{\rho}_t\bar{\theta}_t^{-r} = 0, \qquad \sum_{t=0}^{+\infty}\bar{\rho}_t = +\infty,$$

one can satisfy conditions of the theorem, for example, by taking

$$\rho_t = \max\{\delta_t^\alpha, \bar{\rho}_t\}, \qquad \theta_t = \max\{\delta_t^{\beta/r}, \bar{\theta}_t\},$$

where $0 < \beta < \alpha < 1$. In particular one can take optimal $\alpha, \beta$ from Comment 1.

**Proof of Theorem 2.** Similar to the proof of Theorem 1 we can obtain inequality (18). Taking mathematical expectation from both sides of (18) we get

$$\mathbf{E}(z^{t+1} - g^{t+1})^2 \leq (1 - 2\rho_t)\mathbf{E}(z^t - g^t)^2$$
$$+ 4\rho_t \bar{g} A H \theta_t^\nu + 8\bar{g}\delta_t(x) + 2\rho_t^2 \bar{g}\left(\bar{g} + \bar{K}\theta_t^{-r}\right). \qquad (22)$$

Now the statement follows from Lemma A2 from Appendix.

In the next theorem we assume certain speed of vanishing of the density changes. This can be achieved by increasing frequency of observations $\bar{x}_t$ from densities $g^t(x) = g(\tau_t, x)$ at time moments $\tau_t$ assuming that density $g(\tau, x)$ is Lipschitz continuous with respect to time variable $\tau$. Under these more restrictive assumptions we are able to prove that SQG-estimators $z^t$ track densities $g^t(\cdot)$ with probability one.

**Theorem 3** (*strong consistency of pointwise SQG-estimates*). Assume that

(i) for all $t$ the tracked densities $g^t(x)$ belong to the class $G$ and their changes tend to zero, $|g^{t+1}(x) - g^t(x)| \leq \delta_t(x) \to 0$;

(ii) kernel $K(\cdot)$ is a bounded Borel density satisfying (5), (12);

(iii) sequence of estimates $\{z^t\}$ is obtained by (15), where adjustment coefficients $\rho_t$ and window parameters $\theta_t$ satisfy conditions

$$\lim_{t\to\infty}\frac{\delta_t(x)}{\rho_t}=0, \tag{23}$$

$$\lim_{t\to\infty}\rho_t=\lim_{t\to\infty}\theta_t=0 \text{ a.s.}, \qquad \sum_{t=0}^{+\infty}\rho_t=+\infty \quad \text{a.s.}, \qquad \sum_{t=0}^{+\infty}E\frac{\rho_t^2}{\theta_t^r}<+\infty. \tag{24}$$

Then at any point $x$, where all densities $g^t(x)$ are globally Hölder continuous (see (10)), with probability one holds

$$\lim_{i\to\infty}(z^t-g^t(x))=0.$$

**Comment 3.** Conditions (23) and additional assumption $\delta_t(x)\leq Const\cdot\rho_t$ imply

$$\sum_{t=0}^{+\infty}\delta_t^2\leq \text{Const}\cdot\sum_{t=0}^{+\infty}E\rho_t^2<+\infty$$

and thus assume some rate of convergence $\delta_t\to 0$. So we first have to construct sequences $\rho_t$, $\theta_t$, satisfying (24), choose $\delta_t$ satisfying (23), and then ensure that $|g^{t+1}(x)-g^t(x)|\leq \delta_t$.

**Proof of Theorem 3.** Similar to proof of Theorem 1 we can establish (17). Denote

$$v_t=(z^t-g^t)^2, \quad w_t=2v_t-4\overline{g}\,|g^t-\tilde{g}^t_{\theta_t}(x)|-4\overline{g}\,|g^{t+1}-g^t|-6\overline{g}\,|g^{t+1}-g^t|/\rho_t, \tag{25}$$

$$\gamma_{1t}=\rho_t^2(\xi_t)^2, \quad \gamma_{2t}=2\rho_t\left\langle \xi_t-\nabla\Phi^t_{\theta_t}(z^t),g^{t+1}-z^t\right\rangle.$$

Then (17) can be rewritten in the form

$$v_{t+1}\leq v_t-\rho_t w_t+\gamma_{1t}+\gamma_{2t}.$$

Now let us check conditions of Lemma A3 from Appendix. Condition (i) with $\gamma_t=\gamma_{1t}+\gamma_{2t}$ is satisfied by construction, (ii) is fulfilled by (24). Conditions (iv), (v) are fulfilled by (11), (23) – (25).

Define $\sigma$-algebras $F_t=\sigma\{\overline{x}_0,...,\overline{x}_{t-1}\}$ generated by random variables $\{\overline{x}_0,...,\overline{x}_{t-1}\}$. By (24) the following estimates hold true

$$\mathbf{E}\sum_{t=0}^{+\infty}\gamma_{1t} = \mathbf{E}\sum_{t=0}^{+\infty}\rho_t^2\left(z^t-\frac{1}{\theta_t^r}K\left(\frac{\overline{x}_t-x}{\theta_t}\right)\right)^2$$

$$= 2\mathbf{E}\sum_{t=0}^{+\infty}\rho_t^2\left((z^t)^2+\frac{1}{\theta_t^{2r}}K^2\left(\frac{\overline{x}_t-x}{\theta_t}\right)\right)$$

$$\leq 2\overline{g}^2\sum_{t=0}^{+\infty}\mathbf{E}\rho_t^2+2\overline{K}\sum_{t=0}^{+\infty}\mathbf{E}\rho_t^2\theta_t^{-r}\mathbf{E}\{\theta_t^{-r}K\left(\frac{\overline{x}_t-x}{\theta_t}\right)|F_t\}$$

$$\leq 2\overline{g}^2\sum_{t=0}^{+\infty}\mathbf{E}\rho_t^2+2\overline{g}\overline{K}\sum_{t=0}^{+\infty}\mathbf{E}\rho_t^2\theta_t^{-r}<+\infty,$$

hence $\sum_{i=0}^{+\infty}\gamma_{1i}<+\infty$ a.s.

Random variable $z^t$ is $F_t$-measurable, by (7) $\nabla \Phi_{\theta_t}^t(z^t) = E\{\xi_t | F_t\}$, hence sequence $\{M_t = \sum_{i=0}^{t} \gamma_{2i}, t = 0,1,...\}$ constitutes a martingale with respect to a flow of $\sigma$-algebras $\{F_t\}$. Denote for brevity $\tilde{g}_i = \tilde{g}_{\theta_i}^i(x)$, $K_i = \theta_i^{-r} K((\bar{x}_t - x)\theta_t^{-r})$. By conditions (24),

$$\begin{aligned}
\mathbf{E} M_t^2 &= \sum_{i=0}^{t} \mathbf{E}(\gamma_{2i})^2 \leq 4 \sum_{i=0}^{t} \mathbf{E}\rho_i^2 |g^{t+1} - z^t|^2 \cdot |\tilde{g}_i - K_i|^2 \\
&\leq 8\bar{g}^2 \sum_{i=0}^{t} \mathbf{E}\rho_i^2 (\tilde{g}_i^2 + K_i^2) \leq 8\bar{g}^2 \sum_{i=0}^{t} \mathbf{E}\rho_i^2 (\mathbf{E}\{K_i | F_i\}^2 + K_i^2) \\
&\leq 16\bar{g}^2 \sum_{i=0}^{t} \mathbf{E}\rho_i^2 \mathbf{E}\{K_i^2 | F_i\} \leq 16\bar{g}^2 \bar{K} \sum_{i=0}^{t} \mathbf{E}\rho_i^2 \theta_i^{-r} \mathbf{E}\{K_i | F_i\} \\
&\leq 16\bar{g}^3 \bar{K} \sum_{i=0}^{+\infty} \mathbf{E}\rho_i^2 \theta_i^{-r} < +\infty.
\end{aligned}$$

Hence with probability one, martingale $\{M_t\}$ has a limit, i.e. $\sum_{t=0}^{+\infty} \gamma_{2t} < +\infty$ a.s. This implies that $\sum_{t=0}^{+\infty} \gamma_t = \sum_{t=0}^{+\infty} (\gamma_{1t} + \gamma_{2t}) < +\infty$ a.s., i.e. condition (iii) of Lemma A3 is satisfied. Thus all conditions of Lemma A3 have been verified, that proves the theorem.

### 3.2 Point-wise SQG-estimates in a stationary case

In a stationary situation $g^t(x) \equiv g(x)$ and $|g^{t+1}(x) - g^t(x)| \equiv 0$, but stochastic optimization problems (4) are still nonstationary because involved kernel parameter $\theta$ tends to zero. Convergence results of Theorems 1-3 are applicable to this case, but can be strengthen (prototypes of Theorems 3, 5, 6 for this case were obtained in [30]).

**Theorem 4** (*accuracy of SQG-estimates for small constant adjustment and window parameters*). As before assume that $g(x) \in G$ and kernel $K(\cdot)$ is a bounded Borel density satisfying (5), (12). Let sequence of estimates $\{z^t\}$ is obtained by (15) with constant adjustment and window parameters, $\rho_t \equiv \rho \leq 1/2, \theta_t \equiv \theta > 0$. Then at any point $x$, where density $g(\cdot)$ is globally Hölder continuous (see (10)), the tracking error satisfy

$$\begin{aligned}
(z^{t+1} - g(x))^2 &\leq 2\bar{g}AH\theta^\nu + \bar{g}^2\rho + \bar{K}^2\rho\theta^{-2r} \\
&\quad + (1-2\rho)^t (z^0 - g(x))^2 + \zeta_t,
\end{aligned}$$

where a stochastic term $\zeta_t$ has zero mean $\mathbf{E}\zeta_t = 0$ and a variance $D\zeta_t \leq 2\rho\bar{g}^3(\bar{g} + \bar{K}\theta^{-r})$.

**Comment 4.** Here the error depends on the accuracy of the initial approximation $z^0$, number of steps $t$ and the choice of adjustment parameter $\rho$ and window size $\theta$. The initial error $|z^0 - g(x)|^2$ is quickly suppressed by the multiplier $(1-2\rho)^t$

. For known constants $\bar{g}$, $v$, $H$, $A$, $\bar{K}$ given the number of observations $t$ optimal $\rho$, $\theta$ can be found from minimization of the error function

$$\varphi(\rho,\theta) = 2\bar{g}^2(1-2\rho)^t + 2\bar{g}AH\theta^v + \bar{g}^2\rho + \bar{K}^2\rho\theta^{-2r}$$

over $\rho \in [0,1/2]$, $\theta > 0$. If these constants are not known then a reasonable choice of $\rho$ and $\theta$ can be made assuming that $\rho$, $\theta$ are close to zero. For example, given $\rho$ the asymptotically optimal choice of $\theta$ is $\theta \sim \rho^{\frac{1}{v+2r}}$. Then for large $t$ the deterministic error has the order $\varphi \sim \rho^{\frac{v}{2r+v}}$, and the estimate $\psi$ of the stochastic error (described by the standard deviation $\sqrt{D\zeta_t}$) has the order $\psi \sim \rho^{\frac{r+v}{4(r+v/2)}}$.

Remind that empirical kernel density estimates (13), (14) can be obtained by unconstrained recursive procedure (15) with diminishing adjustment coefficients $\rho_t = 1/t$. So we consider procedures (15) with general rules $\rho_t = \rho/(1+t)^p$, $\theta_t = \theta/(1+t)^q$ with parameters $0 < p \leq 1$, $q \geq 0$. Convergence results of Theorem 2, 3 (with $\delta_t \equiv 0$) hold in this case at any continuity point $x$ of $g(\cdot)$ under weaker than (12) condition (8) due to (9). Moreover we can establish a rate of convergence (in mean) for these procedures and then find optimal parameters $p^*$, $q^*$. It appears that optimal $p^* = 1$ and $0 < q^* < 1$.

**Theorem 5** (*rate of convergence of SQG-estimates for diminishing adjustment and window parameters*). Let $\rho_t = \rho/(1+t)^p$, $\theta_t = \theta/(1+t)^q$, $0 < p \leq 1$, $p < \rho$. Then for any $t \geq 0$ the following estimates hold true

$$\mathbf{E}|z^t - g(x)|^2 \leq Q/(1+t)^\gamma, \quad Q = \max\left\{|z^0 - g(x)|^2, \frac{C}{2\rho - p}\right\}.$$

where $C = \max\{4AH\bar{g}\rho\theta^v, 2\bar{g}^2\rho^2, 2\bar{K}\bar{g}\rho^2\theta^{-r}\}$, $\gamma = \min\{vq, p, (p-rq)\}$.

**Proof of Theorem 5.** Putting $\delta_t \equiv 0$ and $g^t = g^{t+1} = g(x)$ in (22) we have the estimate

$$\mathbf{E}|z^{t+1} - g(x)|^2 \leq (1-2\rho_t)\mathbf{E}|z^t - g(x)|^2$$
$$+ 4AH\bar{g}\rho_t\theta_t^v + 2\rho_t^2\bar{g}^2 + 2\bar{g}\bar{K}\rho_t^2\theta_t^{-r}, \qquad (26)$$

Denoting $v_t = \mathbf{E}\|z^t - g(x)\|^2$ for $\rho_t = \rho/(1+t)^p$, $\theta_t = \theta/(1+t)^q$ for all $t \geq 0$, we get

$$v_{t+1} \leq \left(1 - \frac{2\rho}{(1+t)^p}\right)v_t + \frac{C}{(1+t)^{p+\gamma}},$$

Then by Lemma A3 from Appendix under assumption $\gamma < 2\rho$ we have

$$v_t \leq Q/(1+t)^\gamma, \quad Q = \max\left\{v_0, \frac{C}{2\rho - p}\right\}.$$

**Corollary 1** (*optimal SQG-estimates*). Maximal rate of convergence $1/(1+t)^{\gamma^*}$, $\gamma^* = v/(r+v)$ of SQG-estimates is achieved at $p^* = 1$ and $q^* = 1/(r+v)$. For example, in case $v = r = 1$ maximal rate is proportional to $1/(1+t)^{1/2}$.

### 3.3 Cesàro density estimates in a stationary case

To get more robust estimates we can apply Cesàro averaging to SQG-estimates $z^t$ (15):

$$\bar{z}^t = \sum_{j=0}^{t} z^j \rho_j \bigg/ \sum_{j=0}^{t} \rho_j = (1-\sigma_t)\bar{z}^{t-1} + \sigma_t z^t, \quad \bar{z}^0 = 0, \quad t = 1, 2, \ldots, \quad (27)$$

where $\sigma_t = \rho_t \big/ \sum_{j=0}^{t} \rho_j$. Obviously, if estimates $z^t$ converge to $g(x)$, and $0 \leq \rho_t \leq \rho < \infty$, $\sum_{t=0}^{+\infty} \rho_t = +\infty$, then Cesàro estimates $\bar{z}^t$ also converge to $g(x)$. But for Cesàro estimates we can establish rate of convergence in mean.

**Theorem 6** (*accuracy of Cesàro estimates*). Let density $g(\cdot) \in G$ is Hölder continuous at $x$ (see (11)), $\{\rho_t\}$, $\{\theta_t\}$, be nonnegative deterministic sequences of numbers, constants $\bar{g}(x), v, A, H$ are defined in (2), (11), (12), then for any $t$ holds

$$|\bar{z}^{t+1} - g(x)|^2 \leq \left( \frac{1}{2}(z^0 - g(x))^2 + 2\bar{g}(x)AH \sum_{i=0}^{t} \rho_i \theta_i^v + \bar{g}^2(x) \sum_{i=0}^{t} \rho_i^2 \right.$$

$$\left. + \bar{K}^2 \sum_{i=0}^{t} \rho_i^2 \theta_i^{-2r} \right) \left( \sum_{i=0}^{t} \rho_i \right)^{-1} + \zeta_t, \quad (28)$$

where the stochastic term $\zeta_t$ has zero mean $\mathbf{E}\zeta_t = 0$ and the variance

$$D\zeta_t \leq \bar{g}^3 \left( \sum_{i=0}^{t} \rho_i \right)^{-2} \left( \bar{g} \sum_{i=0}^{t} \rho_i^2 + \bar{K} \sum_{i=0}^{t} \rho_i^2 \theta_i^{-r} \right).$$

**Proof of Theorem 6.** From inequality (18) we obtain:

$$(z^{i+1} - g(x))^2 \leq (z^i - g(x))^2 - 2\rho_i(z^i - g^i)^2$$

$$+ 4\rho_i \bar{g} AH \theta_i^v + 2\rho_i^2 \left( \bar{g}^2 + \bar{K}^2 \theta_i^{-2r} \right) + 2\rho_i \eta_i. \quad (29)$$

where $\eta_i = \left\langle \theta_i^{-r} K(\theta_i^{-1}(\bar{x}_i - x)) - \tilde{g}_{\theta_i}^i(x), g^{i+1} - z^i \right\rangle$. Summing up these inequalities from $i = 0$ to $i = t$ we obtain:

$$0 \leq |z^{t+1} - g(x)|^2 \leq (z^0 - g(x))^2 - 2\sum_{i=0}^{t} \rho_i |z^i - g(x)|^2 +$$

$$+ 4\bar{g}(x)AH \sum_{i=0}^{t} \rho_i \theta_i^v + 2\bar{g}^2(x) \sum_{i=0}^{t} \rho_i^2 + 2\bar{K}^2 \sum_{i=0}^{t} \rho_i^2 \theta_i^{-2r} + \sum_{i=0}^{t} \rho_i \eta_i.$$

Dividing this inequality by $\sum_{i=0}^{t} \rho_i$ and using convexity of function $|z - g(x)|^2$ we obtain the required estimate with $\zeta_t = \sum_{i=0}^{t} \rho_i \eta_i \left( \sum_{i=0}^{t} \rho_i \right)^{-1}$. For the variance we have the estimate

$$D\zeta_t = \left( \sum_{i=0}^{t} \rho_i \right)^{-2} \sum_{i=0}^{t} \rho_i^2 \mathbf{E}\eta_i^2$$

$$\leq \bar{g}^3 \left( \sum_{i=0}^{t} \rho_i \right)^{-2} \left( \bar{g} \sum_{i=0}^{t} \rho_i^2 + \bar{K} \sum_{i=0}^{t} \rho_i^2 \theta_i^{-r} \right).$$

**Corollary 2** (*accuracy of Cesàro estimates for constant adjustment and smoothing parameters*). Let $\theta_i \equiv \theta$, $\rho_i \equiv \rho$, then

$$\left\|\bar{z}^{t+1} - g(x)\right\|^2 \leq \frac{1}{2\rho(t+1)}(z^0 - g(x))^2 + 2\bar{g}(x)AH\theta^v + \bar{g}^2(x)\rho + \bar{K}^2\rho\theta^{-2r} + \zeta_t,$$

where $D\zeta_t \leq \frac{\bar{g}^3}{t+1}\left(\bar{g} + \bar{K}^2\theta^{-r}\right)$.

In Table 1 we compare accuracy of SQG-estimates from Theorem 4 and Cesàro ones from Corollary 2 under the same parameters $\rho$, $\theta$ and at moments $t = \frac{k}{\rho}$, $k = 1, 2, \ldots$. We can see that in the beginning of the estimation process it make sense to apply SQG-estimates since they quickly suppress an initial error, but Cesàro estimates can be applied in the final stage of the estimation process since the lasts suppress a stochastic error.

Table 1. Characteristics of the accuracy of $z^t$, $t = \frac{k}{\rho}$, $k = 1, 2, \ldots$

| Error components | SQG-estimates $z^t$ | Cesàro estimates $z^t$ |
|---|---|---|
| Initial error | $e^{-2k}(z^0 - g(x))^2$ | $\frac{1}{k}(z^0 - g(x))^2$ |
| Deterministic error | $2\bar{g}AH\theta^v + \bar{g}^2\rho + \bar{K}^2\rho\theta^{-2r}$ | $2\bar{g}AH\theta^v + \bar{g}^2\rho + \bar{K}^2\rho\theta^{-2r}$ |
| Stochastic error | $\rho\bar{g}^3(\bar{g} + \bar{K}\theta^{-r})$ | $\frac{1}{k}\bar{g}^3(\bar{g} + \bar{K}\theta^{-r})$ |

### 3.4 General integral SQG-estimates

Point-wise SQG-estimates (15) may be not densities, but they are almost densities on $X$. Indeed, assume for example that $g^t(x) \equiv g(x)$ and density $g(x)$ is continuous for almost all $x \in X \subseteq R^1$. Each estimate $z^t = z^t(x)$ satisfies conditions $0 \leq z^t(x) \leq \bar{g}(x) \leq \bar{g}$. Due to Theorem 3, 5 sequences $\{z^t(x)\}$ almost sure converge to $g(x)$ on everywhere dense subset $X' \subseteq X$, $\lim_t z^t(x) = g(x)$ for all $x \in X$ a.s. Then for any $y$ by Lebesgue dominance convergence theorem

$$\lim_{t\to\infty}\left[G_t(y) = \int_{-\infty}^y z^t(x)dx\right] = \int_{-\infty}^y \lim_{t\to\infty} z^t(x)dx = \int_{-\infty}^y g(x)dx = G(y) \text{ a.s.}$$

To ensure $z^t = z^t(x)$ is a destiny for any $t$, one has include a constraint $\int_X g^t(x)dx = 1$ in the definition of class $G$. The unknown densities may also possess other properties like monotonicity, concavity and etc. It is useful to put corresponding constraints into the estimation problem, the less observation are available the more a priori information should be put into estimation problem to get reasonable estimates (see, discussion in [1, 28, 52]). But in case of general constraints estimation problem becomes not decomposable in $x$.

Instead of (4) consider the following integral estimation problem:

$$\Phi_\theta^t(z(\cdot)) = \frac{1}{2}\int_X\int_X\left[z(x) - \frac{1}{\theta^r}K\left(\frac{\bar{x}-x}{\theta}\right)\right]^2 g^t(\bar{x})h(x)d\bar{x}dx \longrightarrow \min_{z(\cdot)\in G}, \qquad (30)$$

where $h(x)$ is some auxiliary bounded density function on $X$, $\{g^t(\cdot) \in G, t = 0,1,...\}$ are the searched densities, $G$ is a convex set of integrable in square functions on $X$. If $X$ is bounded then density $h$ can be taken uniform. Its $h$-a.s. solution in unconstrained case is:

$$z(x) = \int_X \theta^{-r} K((\bar{x} - x)/\theta) g(\bar{x}) d\bar{x}.$$

In constrained case SQG-estimation procedure is considered in a Hilbert space $L_2(X) = \{z(\cdot) \mid \int_X z^2(x) h(x) dx < +\infty\}$ and has the form:

$$Z^{t+1}(x) = \Pi_G[z^t(x) - \rho_t \xi_t], \quad t = 0,1,..., \tag{31}$$

where stochastic gradient $\xi_t = z^t(x) - \theta_t^{-r} K(\bar{x}_t - x)$, projection operation

$$\Pi_G[z(\cdot)] = argmin_{y(\cdot) \in G} \int_X (z(x) - y(x))^2 h(x) dx, \tag{32}$$

$\{\bar{x}_t\}$ are independent observations from searched densities $g^t(\cdot)$, $t = 0,1,...$.

Consistency of general SQG-estimates is understood in the integral sense:

$$\left\| z^t(\cdot) - g^t(\cdot) \right\|_h^2 = \int_X |z^t(x) - g^t(x)|^2 h(x) dx \to 0 \quad a.s.$$

Due to inequalities

$$\left\| z^{t+1}(\cdot) - g^{t+1}(\cdot) \right\|_h^2 \leq \left\| z^t(\cdot) - g^{t+1}(\cdot) - \rho_t \xi_t \right\|_h^2$$
$$= \int_X |z^t(x) - g^{t+1}(x) - \rho_t \xi_t|^2 h(x) dx$$

one can obtain integral analogs of inequalities (17), (18). Due to possibility to interchange expectation and integration operations one can conclude that (an integral) stochastic term in analog of (18) has zero mean and an integral analog of (20) holds true. Thus all results on convergence and rate of convergence obtained for point-wise case are extended to general estimates (31). The only essential difference is the projection operation in (31). For convergence the set $G$ of feasible densities has to be convex. Similar to [52] objective function and constraints in projection problem (32) can be discretized, then projection problem becomes a finite dimensional quadratic programming problem and can be easily solved.

## 4 Recursive estimation of a vector regression function

The problem of nonparametric regression estimation is to estimate time varying vector regression function $y^t(x)$, $x \in X \subset R^r$, with values in $Y \subset R^m$ from a family of observations on argument $x$ and corresponding output $y$ sampled from joint densities $g^i(x,y), i = 0,1,...,t \to \infty$, i.e. in the present setting we assume that the range $X$ of the argument $x$ is randomly inspected through densities $g_1^t(x) = \int_Y g^t(x,y) dy$. Assume that $y^t(\cdot)$ belong to some class of functions $Y$, in particular it may be $Y = \{y(\cdot) \mid y(x) \in Y(x), x \in X\}$. The searched regression functions $y^t(x)$ are ex-pressed through densities $g^t(x,y)$ as conditional expectations (1). Thus the problem is to track functions $y^t(\cdot)$ through the sequences of independent observations

$$\{(\bar{x}_{01},\bar{y}_{01}),...,(\bar{x}_{0N_0},\bar{y}_{0N_0}), \ ..., \ (\bar{x}_{t1},\bar{y}_{t1}),...,(\bar{x}_{0N_t},\bar{y}_{0N_t})\},$$

and a priori information that $y^t(x)$ belongs to some convex compact set $Y(x) \subseteq R^m$ for any $x \in X$. Each portion of $N_i$ observations $\{(\bar{x}_{t1},\bar{y}_{t1}),...,(\bar{x}_{0N_t},\bar{y}_{0N_t})\}$ is simultaneously and independently made from density $g^i(x,y)$. Beside bounds on $y^t(x)$ the set $Y(x)$ may include relations between components of $y^t(x)$. Assume also that

$$\|Y(x)\| = \sup_{y \in Y(x)} \|y\| < +\infty, \quad g_1(x) \le g_1^t(x) \le \bar{g}_1; \tag{33}$$

$$\int_Y y^2 g_2^t(y)dy \le \bar{g}_3 < +\infty, \quad g_2^t(y) = \int_X g^t(x,y)dx. \tag{34}$$

Consider for $t = 0,1,...$, the following optimization problems:

$$\Phi_\theta^t(x,z) = \frac{1}{2}\int_{X \times Y} \frac{1}{\theta^r} K\left(\frac{\bar{x}-x}{\theta}\right)(z-\bar{y})^2 g^t(\bar{x},\bar{y})d\bar{x}d\bar{y} \to \min_{z \in Y(x)}, \tag{35}$$

and their empirical approximations

$$\Phi_\theta^{N_t}(x,z) = \frac{1}{2N_t}\sum_{i=1}^{N_t} \frac{1}{\theta^r} K\left(\frac{\bar{x}_{ti}-x}{\theta}\right)(z-\bar{y}_{ti})^2 \longrightarrow \min_{z \in Y(x)}, \tag{36}$$

where $\{(\bar{x}_{ti},\bar{y}_{ti})\ i=1,...,N_t\}$ are independently sampled from the same (unknown) density $g^t(x,y)$, $K(\cdot) \ge 0$ is some (bounded and Borel) kernel satisfying (5), (12), and $\theta$ is a positive mollifier (or window) parameter.

Unconstrained minimum in (36) gives a classical Nadaraya-Watson (see [26, 27, 51]) kernel regression estimate

$$\tilde{z}_\theta^{N_t}(x) = \frac{\sum_{i=1}^{N_t} \bar{y}_i K(\frac{\bar{x}_{ti}-x}{\theta})}{\sum_{i=1}^{N_t} K(\frac{\bar{x}_{ti}-x}{\theta})} \tag{37}$$

of the regression function $y^t(x)$.

Beside (35) define functions

$$\Phi^t(x,z) = \frac{1}{2}\int_Y (z-\bar{y})^2 g^t(x,\bar{y})d\bar{y}, \quad t = 0,1,....$$

Unconstrained minimum of $\Phi^t(x,z)$ over $z$ under fixed $x$ is achieved on $y^t(x)$:

$$\Phi^t(x,z) = \frac{1}{2}\int_Y (z-\bar{y})^2 g^t(x,\bar{y})d\bar{y}$$

$$= \frac{1}{2}z^2 \int_Y g^t(x,\bar{y})d\bar{y} - z\int_Y \bar{y}g^t(x,\bar{y})d\bar{y} + \frac{1}{2}\int_Y \bar{y}^2 g^t(x,\bar{y})d\bar{y}$$

$$= \frac{1}{2}g_1^t(x)z^2 - zg_1^t(x)y^t(x) + \frac{1}{2}\int_Y y^2 g^t(x,y)dy$$

$$= \frac{1}{2}g_1^t(x)\left(z-y^t(x)\right)^2 + \frac{1}{2}\int_Y y^2 g^t(x,y)dy - \frac{1}{2}\left(y^t(x)\right)^2 g_1^t(x). \tag{38}$$

Under mild conditions for given $x$ functions $\Phi_\theta^t(x,z)$ uniformly in $y$ converge to $\Phi^t(x,z)$ as $\theta \to 0$.

**Lemma 1.** Assume that $X$ is a compact set, functions $K(\cdot)$ and $g^t(\cdot,\cdot)$ are bounded on $X \times Y$, density $g^t(\cdot, y)$ is continuous at $x$ for almost all $y \in Y$. Then functions $\Phi_\theta^t(x,\cdot)$ uniformly converge on $Y$ to $\Phi^t(x,\cdot)$ as $\theta \to 0$.

**Proof.** Changing variables in (35) we obtain

$$\Phi_\theta^t(x,z) = \frac{1}{2}\int_{R^r \times Y}(z-\bar{y})^2 K(\tilde{x}) g^t(x+\theta\tilde{x},\bar{y})d\tilde{x}d\bar{y}.$$

Then by Lebesgue dominance convergence theorem

$$\lim_{\theta \to 0}\Phi_\theta^t(x,z) = \frac{1}{2}\int_Y (z-\bar{y})^2 \lim_{\theta \to 0}\left(\int_{R^r} K(\tilde{x})g^t(x+\theta\tilde{x},\bar{y})d\tilde{x}\right)d\bar{y}$$

$$= \frac{1}{2}\int_Y (z-\bar{y})^2 g^t(x,\bar{y})\int_X K(\tilde{x})d\tilde{x}d\bar{y}$$

$$= \frac{1}{2}\int_Y (z-\bar{y})^2 g^t(x,\bar{y})d\bar{y} = \Phi^t(x,z).$$

Since $\Phi_\theta^t(x,\cdot)$ is convex then a uniform convergence follows from the pointwise one.

**Lemma 2.** Assume that for almost all $y \in Y$ densities $g^t(x,y)$ are globally Hölder continuous in $x$:

$$|g^t(x_1,y) - g^t(x_2,y)| \le L(y)\|x_1 - x_2\|^\nu, \tag{39}$$

with constants $\nu$, $L(y)$ such that (8) holds and

$$\int_Y L(y)dy = B < +\infty, \quad \int_Y \|y\|^2 L(y)dy = C < +\infty. \tag{40}$$

Then the following estimates hold

$$|\Phi^t(x,z) - \Phi_\theta^t(x,z)| \le A(\|z\|^2 B + C)\theta^\nu, \tag{41}$$

**Proof.** Indeed,

$$|\Phi^t(x,z) - \Phi_\theta^t(x,z)| \le \frac{1}{2}\int_{X \times Y}(z-\bar{y})^2 K(\tilde{x})|g^t(x+\theta\tilde{x},\bar{y}) - g^t(x,\bar{y})|d\tilde{x}d\bar{y}$$

$$\le \frac{1}{2}\theta^\nu \int_X \|x\|^\nu K(x)dx \cdot \int_Y (z-\bar{y})^2 L(\bar{y})d\bar{y}$$

$$\le A(\|z\|^2 B + C)\theta^\nu.$$

## 4.1 Point-wise SQG-estimates of a vector regression function

In the present paper we consider constrained kernel regression estimates obtained from a family of time varying stochastic programming problems (35) with $\theta_t \to 0$, $t \to +\infty$ (the so-called nonstationary stochastic optimization problem). Problems (35) are not given explicitly because of unknown density $g^t(\cdot,\cdot)$. Empirical approximations (36) also may not be accurate because of insufficient number $N_t$ of observations made from the same density $g^t(\cdot,\cdot)$. So we construct recursive estimates, namely at iteration $t$ we update previously obtained regression estimate $z^t(\cdot)$ by one step of minimization of (36). Remark that problems (36) have to be considered for many $x$ to get a table for the function values $y^t(x)$. For more or less complex constraints $Y(x)$ this may be a tedious task. From this point of view relatively simple recursive estimates also may be helpful. Thus we solve nonstationary problems (35) with $\theta \to 0$ by iterative SQG-method [6, 7, 8] and obtain a sequence of SQG-estimates:

$$z^{t+1} = \Pi_{Y(x)}(z^t - \rho_t \xi_t), \quad z^0 \in Y(x), \quad i = 0, 1, \ldots, \tag{42}$$

with stochastic gradients

$$\xi_t = \frac{1}{\theta_t^r N_t} \sum_{i=1}^{N_t} K\left(\frac{\overline{x}_{ti} - x}{\theta_t}\right) \cdot (z^t - \overline{y}_{ti}) \tag{43}$$

such that conditional expectation $E\{\xi_t \mid z^t\} = \nabla \Phi_{\theta_t}(x, z^t)$, nonnegative smoothing parameters $\theta_t$ and adjustment coefficients $\rho_t$ are measurable with respect to $\sigma$-algebra $\sigma_{t-1}$ generated by i.i.d. observations $\{(\overline{x}_{01}, \overline{y}_{01}), \ldots, (\overline{x}_{t-1,N_{t-1}}, \overline{y}_{t-1,N_{t-1}})\}$ defined on a joint probability space $(\Omega, \Sigma, P)$. Here estimates $z^t$ and stochastic gradients $\xi^t$ depend on $x$ but for brevity we drop the argument.

For the stationary case $g^t(\cdot,\cdot) \equiv g(\cdot,\cdot)$ similar recursive estimation procedure but without projection on the set $Y(x)$ and with deterministic sequences $\{\theta_t, \rho_t\}$ were considered in [44, 46]. Remark that constrained estimates $z^t$ of the vector regression function $y^t(x)$ are nonlinear in observations $\{\overline{y}_0, \overline{y}_1, \ldots, \overline{y}_t\}$.

Next theorems establish rate of tracking $y^t(x)$ by estimators $z^t$ under different assumptions on the speed of regression changes, and hence for the stationary case.

**Theorem 7** (*accuracy of regression estimates*). Assume that

(i) for all $t$ the tracked regressions $y^t(x) \in Y(x)$ and their changes satisfy $\|y^{t+1}(x) - y^t(x)\| \leq \delta(x)$;

(ii) all densities $g^t(x, y)$ are globally Hölder continuous in $x$ (see (39)) and satisfy conditions (33), (34);

(iii) for a given $x$ all densities $g^t(x, \cdot)$ have bounded supports, $\sup\{\|y\| \mid g^t(x, y) > 0\} \leq \|Y_x\| < +\infty$;

(iv) kernel $K(\cdot)$ is a bounded Borel density satisfying (5), (12);

(v) sequence of estimates $\{z^t\}$ is obtained by (42), where adjustment coefficients $\rho_t$ and window parameters $\theta_t$ are constant, $\rho_t \equiv \rho \leq 1/\overline{g}_1, \theta_t \equiv \theta$.

Then at any point $x$, where $g_1^t(x) \geq g_1(x) > 0$, the tracking error satisfy

$$(z^{t+1} - y^{t+1}(x))^2 \leq (1 - g_1(x)\rho)^t \cdot (z^0 - y^0(x))^2$$
$$+ Q(x) g_1^{-1}(x)\left(\delta \rho^{-1} + \delta + \theta^\gamma + \rho \theta^{-2r}\right) + 2\zeta_t, \tag{44}$$

where a stochastic term $\zeta_t$ has zero mean $\mathbf{E}\zeta_t = 0$ and a variance
$$D\zeta_t \leq \rho\theta^{-r}P(x),$$
$P(x)$ and $Q(x)$ are some constants (defined in (49), (54)).

**Comment 5** (*error minimization*). The regression error estimate is quite similar to the one for density estimates (see Theorem 1). The error depends on the accuracy of the initial approximation $z^0$, number of steps $t$, the speed of changes $\delta$ and the choice of adjustment parameter $\rho$ and window size $\theta$. So an asymptotically optimal choice of $\rho$ and $\theta$ can be made as in Comment 1. The main difference is the presence of factor $g_1(x)$ in (44). This may cause additional requirement on the choice of $\rho < 1/\overline{g}_1$ and worsen the estimates.

**Proof of Theorem 7.** Denote for brevity $y^t = y^t(x)$. Below $\langle \cdot, \cdot \rangle$ designates an inner product of two vectors. Estimate

$$\begin{aligned}
(z^{t+1} - y^{t+1})^2 &\leq (z^t - y^{t+1} - \rho_t\xi_t)^2 \\
&\leq (z^t - y^{t+1})^2 - 2\rho_t\langle \xi_t, z^t - y^{t+1}\rangle + \rho_t^2(\xi_t)^2 \\
&\leq (z^t - y^t)^2 + 2\langle z^t - y^t, y^t - y^{t+1}\rangle + (y^t - y^{t+1})^2 + \rho_t^2(\xi_t)^2 \\
&\quad -2\rho_t\langle \nabla_z\Phi_{\theta_t}^t(x,z^t), z^t - y^{t+1}\rangle \\
&\quad -2\rho_t\langle \xi_t - \nabla_z\Phi_{\theta_t}^t(x,z^t), z^t - y^{t+1}\rangle \\
&\leq (z^t - y^t)^2 + 2\langle z^t - y^t, y^t - y^{t+1}\rangle + (y^t - y^{t+1})^2 \\
&\quad -2\rho_t\langle \nabla_z\Phi_{\theta_t}^t(x,z^t), z^t - y^t\rangle - 2\rho_t\langle \nabla_z\Phi_{\theta_t}^t(x,z^t), y^t - y^{t+1}\rangle \\
&\quad -2\rho_t\langle \xi_t - \nabla_z\Phi_{\theta_t}^t(x,z^t), z^t - y^{t+1}\rangle + \rho_t^2(\xi_t)^2
\end{aligned}$$

From (38) follows $\Phi^t(x,z) - \Phi^t(x,y^t) = \frac{1}{2}g_1^t(x)(z - y^t)^2$. Denote
$$\Delta_t(x,\theta) = \sup_{y \in Y(x)} |\Phi^t(x,y) - \Phi_\theta^t(x,y)|.$$

Then by convexity of $\Phi_{\theta_t}^t(x,\cdot)$ we have

$$\begin{aligned}
-\langle \nabla_z\Phi_{\theta_t}^t(x,z^t), z^t - y^t\rangle &\leq \Phi_{\theta_t}^t(x,y^t) - \Phi_{\theta_t}^t(x,z^t) \\
&\leq \Phi^t(x,y^t) - \Phi^t(x,z^t) \\
&\quad + 2\sup_{y \in Y(x)} |\Phi^t(x,y) - \Phi_{\theta_t}^t(x,y)| \\
&= -\frac{1}{2}g_1^t(x)(z^t - y^t)^2 + 2\Delta_t(x,\theta_t).
\end{aligned}$$

So,
$$\begin{aligned}
(z^{t+1} - y^{t+1})^2 &\leq (z^t - y^t)^2 - \rho_t g_1^t(x)(z^t - y^t)^2 + 2\langle z^t - y^t, y^t - y^{t+1}\rangle \\
&\quad + (y^t - y^{t+1})^2 - 2\rho_t\langle \nabla_z\Phi_{\theta_t}^t(x,z^t), y^t - y^{t+1}\rangle + 2\Delta_t(x,\theta_t) \\
&\quad + \rho_t^2(\xi_t)^2 - 2\rho_t\langle \xi_t - \nabla_z\Phi_{\theta_t}^t(x,z^t), z^t - y^{t+1}\rangle.
\end{aligned}$$

Taking into account bounds (33), the expression $\nabla_z \Phi^t(x,z) = g_1^t(x)(z-y^t)$ and the result of Lemma 2, from the above inequality we get:

$$(z^{t+1} - y^{t+1})^2 \leq (z^t - y^t)^2 - g_1^t(x)\rho_t(z^t - y^t)^2 + 2\|z^t - y^t\| \cdot \|y^t - y^{t+1}\|$$
$$+ \|y^t - y^{t+1}\|^2 + 2\rho_t \|\nabla_z \Phi_{\theta_t}^t(x,z^t)\| \cdot \|y^t - y^{t+1}\|$$
$$+ 4\rho_t \Delta_t(x,\theta_t) + \rho_t^2(\xi_t)^2 - 2\rho_t \langle \xi_t - \nabla_z \Phi_{\theta_t}^t(x,z^t), z^t - y^{t+1}\rangle$$
$$\leq (z^t - y^t)^2 - g_1(x)\rho_t(z^t - y^t)^2$$
$$+ 6\|Y(x)\|\delta_t(x) + 4\rho_t\|Y(x)\|\bar{g}_1\delta_t(x) + \rho_t^2(\xi_t)^2$$
$$+ 4\rho_t \Delta_t(x,\theta_t) + 2\rho_t \eta_t, \quad (45)$$
$$\leq (1 - g_1(x)\rho_t)(z^t - y^t)^2$$
$$+ 6\|Y(x)\|\delta_t(x) + 4\rho_t\|Y(x)\|\bar{g}_1\delta_t(x) + \rho_t^2(\xi_t)^2$$
$$+ 4\rho_t \theta_t^\nu A\left(B\|Y(x)\|^2 + C\right) + 2\rho_t \eta_t, \quad (46)$$

where $\eta_t = \langle \nabla_z \Phi_{\theta_t}^t(x,z^t) - \xi_t, z^t - y^{t+1}\rangle$. Estimate

$$(\xi_t)^2 \leq \frac{2}{\theta_t^{2r}}\bar{K}^2\left(\|z^t\|^2 + \frac{1}{N_t}\sum_{i=1}^{N_t}\|\bar{y}_{ti}\|^2\right) \leq \frac{2}{\theta_t^{2r}}\bar{K}^2(\|Y(x)\|^2 + \|Y_x\|^2). \quad (47)$$

Using (47) from (46) we obtain
$$(z^{t+1} - y^{t+1})^2 \leq (1 - g_1(x)\rho_t)(z^t - y^t)^2 +$$
$$+ 6\|Y(x)\|\delta_t(x) + 4\rho_t\|Y(x)\|\bar{g}_1\delta_t(x)$$
$$+ 4\rho_t\theta_t^\nu A\left(B\|Y(x)\|^2 + C\right) + 2\rho_t^2\theta_t^{-2r}\bar{K}^2\left(\|Y(x)\|^2 + \|Y_x\|^2\right)$$
$$+ 2\rho_t \eta_t. \quad (48)$$

Denote $v_t = (z^t - y^t)^2$,
$$Q(x) = \max\left\{6\|Y(x)\|, 4\|Y(x)\|\bar{g}, 4A\left(B\|Y(x)\|^2 + C\right), 2\bar{K}^2\left(\|Y(x)\|^2 + \|Y_x\|^2\right)\right\}, \quad (49)$$

and assume $\rho_t \equiv \rho$, $\theta_t \equiv \theta$, $\delta_t(x) \leq \delta$, then (48) becomes
$$v_{t+1} \leq (1 - g_1(x)\rho)v_t + Q(x)\left(\delta + \rho\delta + \rho\theta^\nu + \rho^2\theta^{-2r}\right) + 2\rho_t \eta_t.$$

For $\rho < 1/g_1(x)$ by Lemma A1 from Appendix we get
$$v_{t+1} \leq (1 - g_1(x)\rho)^t v_0 + Q(x)g_1^{-1}(x)\left(\delta\rho^{-1} + \delta + \theta^\nu + \rho\theta^{-2r}\right) + 2\zeta_t. \quad (50)$$

where a stochastic term $\zeta_t = \sum_{i=0}^{t}\rho(1 - 2g_1(x)\rho)^{t-i}\eta_i$ has zero mean $\mathbf{E}\zeta_t = 0$ and the variance

$$D\zeta_t = \sum_{i=0}^{t}\rho^2(1 - 2g_1(x)\rho)^{2(t-i)}E\eta_i^2 \leq \frac{\rho}{2g_1(x)}\sup_i E\eta_i^2,$$

Denote $\psi_{ij} = \frac{1}{\theta_i^r}K\left(\frac{\bar{x}_{ij} - x}{\theta_i}\right) \cdot (z^i - \bar{y}_{ij}))$, and estimate

$$\mathbf{E}\{\|\psi_{ij}\|^2 \mid z^i\} = \int_{X \times Y}(z^i - \bar{y})^2 \theta_i^{-2r}K^2((\bar{x} - x)\theta_i^{-1})g^i(\bar{x}, \bar{y})d\bar{x}d\bar{y}$$

$$\leq \theta_i^{-r} \bar{K} \int_{X \times Y} (z^i - \bar{y})^2 \theta_i^{-r} K((\bar{x} - x)\theta_i^{-1}) g^i(\bar{x}, \bar{y}) d\bar{x} d\bar{y}$$

$$\leq 2\theta_i^{-r} \bar{K} \left( \|z^i\|^2 \int_{X \times Y} \theta_i^{-r} K((\bar{x} - x)\theta_i^{-1}) g^i(\bar{x}, \bar{y}) d\bar{x} d\bar{y} \right.$$

$$\left. + \int_{X \times Y} \theta_i^{-r} K((\bar{x} - x)\theta_i^{-1}) \bar{y}^2 g^i(\bar{x}, \bar{y}) d\bar{x} d\bar{y} \right)$$

$$\leq 2\bar{K} \left( \sup_{y \in Y(x)} \|y\|^2 \sup_{x \in X} g_1^i(x) + \sup_{x \in X} \int_Y \bar{y}^2 g^i(x, \bar{y}) d\bar{y} \right) \theta_i^{-r}$$

$$\leq 2\bar{K} \left( \|Y(x)\|^2 \bar{g}_1 + \bar{g}_3 \right) \theta_i^{-r}. \tag{51}$$

Then

$$\mathbf{E}\eta_i^2 = \mathbf{E}\left\langle z^i - y^i, \nabla_z \Phi_{\theta_i}^i(x, z^i) - \xi_i \right\rangle^2 = \mathbf{E}\left\langle z^i - y^i, \nabla_z \Phi_{\theta_i}^i(x, z^i) - \frac{1}{N_i} \sum_{j=1}^{N_i} \psi_{ij} \right\rangle^2$$

$$= \mathbf{E}\left( \frac{1}{N_i} \sum_{j=1}^{N_i} \left\langle z^i - y^i, \nabla_z \Phi_{\theta_i}^i(x, z^i) - \psi_{ij} \right\rangle \right)^2 = \frac{1}{N_i^2} \sum_{j=1}^{N_i} \mathbf{E}\left\langle z^i - y^i, \nabla_z \Phi_{\theta_i}^i(x, z^i) - \psi_{ij} \right\rangle^2$$

$$= \frac{1}{N_i} \mathbf{E}\left\langle z^i - y^i, \nabla_z \Phi_{\theta_i}^i(x, z^i) - \psi_{i1} \right\rangle^2 \leq \frac{4}{N_i} \mathbf{E}(\|z^i\|^2 + \|y^i\|^2) \cdot \left( \|\nabla_z \Phi_{\theta_i}^i(x, z^i)\|^2 + \psi_{i1}^2 \right)$$

$$\leq \frac{16\|Y(x)\|^2}{N_i} \mathbf{E}\psi_{i1}^2 = \frac{16\|Y(x)\|^2}{N_i} \mathbf{E}\mathbf{E}\{\psi_{i1}^2 \mid z^i\} = \frac{32\|Y(x)\|^2}{N_i \theta_i^r} \bar{K}\left( \|Y(x)\|^2 \bar{g}_1 + \bar{g}_3 \right), \tag{52}$$

and

$$\mathbf{E}\xi_t^2 = \mathbf{E}\left[ \frac{1}{N_t} \sum_{j=1}^{N_t} \psi_{tj} \right]^2 = \frac{1}{N_t^2} \sum_{j=1}^{N_t} \mathbf{E}\psi_{tj}^2 = \frac{1}{N_t} \mathbf{E}\psi_{t1}^2 \leq \frac{2}{N_t \theta_t^r} \bar{K}\left( \|Y(x)\|^2 \bar{g}_1 + \bar{g}_3 \right). \tag{53}$$

Thus

$$D\zeta_t \leq \rho \theta^{-r} \frac{16\bar{K}\|Y(x)\|^2 (\|Y(x)\|^2 \bar{g}_1 + \bar{g}_3)}{g_1(x) \inf_{i \leq t} N_i} = \rho \theta^{-r} P(x). \tag{54}$$

In the next theorem we assume that regression changes vanish with time.

**Theorem 8** (*convergence in mean of regression estimates*). Assume that

(i) for all $t$ the tracked regressions $y^t(x) \in Y(x)$ and their changes vanish with rate $\delta_t(x) \to 0$, i.e. $\|y^{t+1}(x) - y^t(x)\| \leq \delta_t(x) \to 0$ as $t \to \infty$;

(ii) all densities $g^t(x, y)$ are globally Hölder continuous in $x$ (see (39)) and satisfy conditions (33), (34);

(iii) kernel $K(\cdot)$ is a bounded Borel density satisfying (5), (12);

(iv) sequence of estimates $\{z^t\}$ is obtained by (42), where adjustment coefficients $\rho_t$ and window parameters $\theta_t$ are deterministic and satisfy conditions (21).

Then at any point $x$, where all $g_1^t(x) \geq g_1(x) > 0$, holds

$$\lim_{t \to \infty} \mathbf{E}(z^t - y^t(x))^2 = 0.$$

**Proof.** Similar to proof of Theorem 7 we can obtain inequality (45). Taking expectations from both sides of (45) and accounting for (52), we obtain

$$\mathbf{E}(z^{t+1} - y^{t+1})^2 \leq (1 - g_1(x)\rho_t)\mathbf{E}(z^t - y^t)^2 + 6\|Y(x)\|\delta_t(x)$$
$$+ 4\rho_t\|Y(x)\|\bar{g}_1\delta_t(x) + 4\rho_t\theta_t^\gamma A\left(B\|Y(x)\|^2 + C\right)$$
$$+ \frac{2\rho_t^2}{N_t\theta_t^r}\bar{K}\left(\|Y(x)\|^2\bar{g}_1 + \bar{g}_3\right). \tag{55}$$

Now the statement $v_t = \mathbf{E}(z^t - y^t)$ follows from Lemma A2 from Appendix.

To prove strong consistency of regression estimates we have to assume existence of certain rate of nonstationarity $\delta_t(x) \to 0$.

**Theorem 9** (strong pointwise consistency of recursive regression estimates). In conditions of Theorem 8 assume that instead of (21) $\rho_t$, $\theta_t$ and $\delta_t$ satisfy (23), (24). Then at any point $x$, where all $g_1^t(x) \geq g_1(x) > 0$, with probability one holds

$$\lim_{t \to \infty}(z^t - y^t(x)) = 0.$$

The proof is similar to the one of Theorem 3 and is based on inequalities (45), (52), (53) and Lemma A4 from Appendix.

**Comment 6.** Theorems 7 - 9 are valid for the stationary case when $g^t(x, y) \equiv g(x, y)$ and hence $y^t(x) \equiv y(x)$. Moreover convergence results of Theorems 8, 9 due to Lemma 1 hold true under weaker than (39), (40) assumption that joint stationary density $g(x, y)$ is continuous at $x$ for almost all $y$. In such a case density $g_1(x)$ and regression function $y(x)$ are continuous at $x$.

### 4.2 Cesàro estimates of a stationary regression function

In a stationary situation $g^t(x, y) \equiv g(x, y)$ and $y^t(x) \equiv y(x)$. In this case we can apply Cesàro averaging to SQG-estimates (42):

$$\bar{z}^t = \sum_{i=0}^{t}\rho_i z^i \bigg/ \sum_{i=0}^{t}\rho_i = (1 - \sigma_t)\bar{z}^{t-1} + \sigma_t z^t, \quad \bar{z}^0 = 0, \quad t = 0, 1, ..., \tag{56}$$

where $\sigma_t = \rho_t \big/ \sum_{i=0}^{t}\rho_i$. Obviously, if estimates $z^t$ converge to $y(x)$ a.s., and $0 \leq \rho_t \leq \rho < \infty$, $\sum_{t=0}^{\infty}\rho_t = +\infty$, then Cesàro estimates $\bar{z}^i$ also converge to $y(x)$ a.s. For Cesàro estimates we can establish convergence and rate of convergence.

**Lemma 3** (*accuracy of Cesàro estimates*). Let $\{\theta_i\}$, $\{\rho_i\}$ be deterministic numbers Then for any $t$ holds

$$g_1(x)\|\bar{z}^t - y(x)\|^2 \leq \left((z^0 - y^0)^2 + 2\sum_{i=0}^{t}\rho_i\Delta_i(x)\right.$$
$$\left. + \sum_{i=0}^{t}\rho_i^2\|\xi_i\|^2 + 2\sum_{i=0}^{t}\rho_i\eta_i\right) \cdot \left(\sum_{i=0}^{t}\rho_i\right)^{-1}. \tag{57}$$

where

$$\Delta_i(x) = \sup_{y \in Y(x)} |\Phi_{\theta_i}(x,y) - \Phi(x,y)|,$$

$$(\xi_i)^2 \le \frac{2}{\theta_i^{2r}} \bar{K}^2 (\|Y(x)\|^2 + \|Y_x\|^2),$$

$$\mathbf{E}\xi_i^2 \le \frac{2}{N_i \theta_i^r} \bar{K}\left(\|Y(x)\|^2 \bar{g}_1 + \bar{g}_3\right),$$

$$\mathbf{E}\eta_i = 0, \quad E\eta_i^2 \le \frac{32\|Y(x)\|^2}{N_i \theta_i^r} \bar{K}\left(\|Y(x)\|^2 \bar{g}_1 + \bar{g}_3\right).$$

**Proof.** From (45) with $\delta_t(x) \equiv 0$ we have

$$(z^{i+1} - y^{i+1})^2 \le (1 - g_1(x)\rho_t)\cdot(z^i - y^i)^2 + 2\rho_i\Delta_i(x,\theta_i) + \rho_i^2(\xi_i)^2 + 2\rho_i\eta_i,$$

where $\eta_i = \langle \nabla_z \Phi_{\theta_i}^i(x,z^t) - \xi_i, z^i - y^{i+1}\rangle$. Summing up these inequalities from $i=0$ to $i=t$, dividing the result by $\sum_{i=0}^t \rho_i$ and using Jensen inequality

$$\|\bar{z}^t - y(x)\|^2 \le \sum_{i=1}^t \rho_i \|y^i - y(x)\|^2 / \sum_{i=0}^t \rho_i,$$

we obtain the required estimate (57). Estimates for $(\xi_i)^2$, $\mathbf{E}(\xi_i)^2$, $\mathbf{E}(\eta_i)^2$ were obtained in (47), (53), (52) respectively.

**Corollary 3** (*convergence in mean of Cesáro estimates*). Assume that density $g(\cdot,y)$ is continuous at $x$ for almost all $y$, parameters $\rho_t$, $\theta_t$ satisfy (21). Then Cesáro estimates $\bar{z}^t$ converge to $y(x)$ in mean, i.e. $\lim_{t\to\infty} g_1(x)E\|\bar{z}^t - y(x)\|^2 = 0$.

The statement follows from estimate (57), Lemma 1 and (21).

**Corollary 4** (*accuracy of Cesáro estimates for constant adjustment and smoothing parameters*). Let $\rho_i \equiv \rho$, $\theta_i \equiv \theta$, and suppose that density $g(x,y)$ is globally Hölder continuous at $x$ (see (39), (40)), then

$$g_1(x)\mathbf{E}\|\bar{z}^{t+1} - y(x)\|^2 \le \left(\frac{(z^0 - y(x))}{\rho t} + Q_1(x)\theta^\nu + Q_2(x)\frac{\rho}{\theta^{-r}\min_{i\le t} N_i}\right),$$

$$g_1(x)\|\bar{y}^i - y(x)\|^2 \le \left(\frac{(z^0 - y(x))}{\rho t} + Q_1(x)\theta^\nu + Q_3(x)\rho\theta^{-2r}\right) + \zeta_t,$$

where factors $Q_1(x) = 4A\left(B\|Y(x)\|^2 + C\right)$, $Q_2(x) = 2\bar{K}\left(\|Y(x)\|^2 \bar{g}_1 + \bar{g}_3\right)$, $Q_3(x) = 2\bar{K}^2(\|Y(x)\|^2 + \|Y_x\|^2)$, stochastic tern $\zeta_t$ has zero mean and a variance $D\zeta_t \le 32\|Y(x)\|^2 Q_2(x)\frac{\rho}{t\theta^r \min_{i\le t} N_i}$.

**Comment 7.** Similar to density estimates (see Table 1) accuracy of SQG-estimates (Theorem 7) and Cesáro ones (Corollary 4) of regression function $y(x)$ have similar deterministic component, but differently diminishing initial error and stochastic components.

General regression SQG-estimates from some class $Y$ reflecting properties of regression $y^t(x)$ can be constructed similar to subsection 3.4.

## 5 Conclusions

We constructed constrained recursive SQG-estimators which track possibly moving density and regression functions. These estimators are of Monte Carlo type, i.e. they sequentially use the samples from unknown changing densities and random input/output regression data. We evaluated accuracy of the obtained estimates at any iteration and for any adjustment (step) and kernel (mollifier) parameters, including not tending to zero. Accuracy of estimates are three fold, first they account for inaccuracy of initial approximation, second they exhibit deterministic systematic bias, and third describe stochastic errors. The obtained accuracy estimates allowed finding asymptotically optimal adjustment and kernel parameters of the estimation procedures. In case of known bound on the speed of nonstationarity adjustment and kernel parameters are adapted to this speed. Asymptotically the resulting error is expressed in terms of the speed of nonstationarity. In a stationary case the kernel parameter is adapted to the adjustment parameter, and the latter is selected from a tradeoff between the asymptotic systematic and disappearing initial error. In a stationary case we also consider the so-called Cesáro estimators obtained from SQG-estimates by means of Cesáro averaging. Because of different accuracy characteristics they can be applied at the later stages of the estimation process.

Informally speaking the results of the paper show that unknown even changing objects (in our case distributions and dependencies) can be perceived (to certain extent) through random observations over them.

To construct and investigate SQG-estimates we employed some results from the field of nonstationary constrained stochastic optimization and stochastic qusi-gradient (SQG) methods. Other ideas from these fields, like adaptation of parameters, multi-step optimization methods, may be also applied to improve estimation.

## 6 Appendix

This section contains useful lemmas on recursively estimated number sequences, convergent toward zero.

**Lemma A1.** Assume that number sequence $\{v_t\}$ admits the following reqursive estimation
$$v_{t+1} \leq (1 - q\rho_t)v_t + \rho_t \alpha_t + \rho_t \eta_t, \quad t \geq t_0,$$
where $q \geq 0$, $q\rho_t \leq 1$ and $0 \leq \alpha_t \leq \alpha$. Then

$$v_{t+1} \leq v_{t_0}\prod_{i=t_0}^{t}(1-q\rho_i)+\frac{\alpha}{q}+\sum_{i=t_0}^{t}\eta_i\rho_i\prod_{j=i+1}^{t}(1-q\rho_j)$$

$$\leq v_{t_0}\exp\{-q\sum_{i=t_0}^{t}\rho_i\}+\frac{\alpha}{q}+\sum_{i=t_0}^{t}\eta_i\rho_i\prod_{j=i+1}^{t}(1-q\rho_j).$$

In particular, if

$$v_{t+1}\leq(1-q\rho)v_t+\rho(\alpha+\eta_t),\quad t\geq t_0,$$

then

$$v_{t+1}\leq v_{t_0}(1-q\rho)^{t-t_0}+\frac{\alpha}{q}+\sum_{i=t_0}^{t}\rho(1-q\rho)^{t-i}\eta_i.$$

**Proof.** The estimates follow from inequalities:

$$v_{t+1} \leq v_{t_0}\prod_{i=t_0}^{t}(1-q\rho_i)+\sum_{i=t_0}^{t}(\alpha+\eta_i)\rho_i\prod_{j=i+1}^{t}(1-q\rho_j)$$

$$= v_{t_0}\prod_{i=t_0}^{t}(1-q\rho_i)+\sum_{i=t_0}^{t}\eta_i\rho_i\prod_{j=i+1}^{t}(1-q\rho_j)$$

$$+\frac{\alpha}{q}\sum_{i=t_0}^{t}\left(\prod_{j=i+1}^{t}(1-q\rho_j)-\prod_{j=i}^{t}(1-q\rho_j)\right)$$

$$\leq v_{t_0}\prod_{i=t_0}^{t}(1-q\rho_i)+\sum_{i=t_0}^{t}\eta_i\rho_i\prod_{j=i+1}^{t}(1-q\rho_j)$$

$$+\frac{\alpha}{q}\left(1-\prod_{j=t_0}^{t}(1-q\rho_j)\right)$$

$$\leq v_{t_0}\prod_{i=t_0}^{t}(1-q\rho_i)+\frac{\alpha}{q}+\sum_{i=t_0}^{t}\eta_i\rho_i\prod_{j=i+1}^{t}(1-q\rho_j)$$

$$\leq v_{t_0}\exp\{-q\sum_{i=t_0}^{t}\rho_i\}+\frac{\alpha}{q}+\sum_{i=t_0}^{t}\eta_i\rho_i\prod_{j=i+1}^{t}(1-q\rho_j)$$

The following lemma strengthens the result of Lemma A1 under additional assumption on $\alpha_t$.

**Lemma A2** [32]. Let

$$v_{t+1}\leq(1-\rho_t)v_t+\rho_t\alpha_t,\quad 0<\rho_t\leq 1,\quad \alpha_t\geq 0;$$

$$\sum_{t=0}^{\infty}\rho_t=+\infty,\quad \lim_{t\to\infty}\alpha_t=0.$$

Then $\overline{\lim}_{t\to\infty}v_t\leq 0$. In particular if $v_t\geq 0$ then $\lim_{t\to\infty}v_t=0$.

In conditions of Lemma A2 one can estimate rate of convergence of $v_t$ to zero.

The following lemma is an analog of well known Chung's lemma (see [32]), but with estimate valid for all $t$.

**Lemma A3** [19, p.282]. Let $\{v_i\}_{i\geq 0}$ be a sequence of nonnegative numbers such that

$$v_{t+1} \leq \left(1 - \frac{\rho}{(1+t)^\beta}\right) v_t + \frac{C}{(1+t)^{\beta+\gamma}}, \quad v_0 < +\infty,$$

where $0 < \beta \leq 1$, $0 < \gamma < \rho$, $C > 0$. Then for any $t \geq 0$

$$v_t \leq \frac{Q}{(1+t)^\gamma}, \quad Q = \max\left\{v_0, \frac{C}{\rho - \beta}\right\}.$$

We also use some stochastic version of Lyapunov function method for discrete time dynamic processes (Lemma A3). Its essential feature is that Lyapunov function values $v_t$ do not necessarily decrease along the trajectory of the process, i.e. corresponding variable $w_t$ below is not necessarily nonnegative.

**Lemma A4** [10]. Let $v_t \geq 0$, $\rho_t \geq 0$, $w_t$, $\gamma_t$, $t \geq 0$, be a sequence of random variables (scalars). Suppose each of the following conditions is fulfilled with probability one:

(i) $v_{t+1} \leq v_t - \rho_t w_t + \gamma_t$, all $t \geq 0$;

(ii) $\lim_{t \to \infty} \rho_t = 0$, $\sum_{t=0}^{\infty} \rho_t = +\infty$;

(iii) $v_t < +\infty$ all $t \geq 0$, $\sum_{t=0}^{\infty} \gamma_t < +\infty$;

(iv) for any indices $\{t_s \to \infty\}$ if $\liminf_s v_{t_s} > 0$ then $\liminf_s w_{t_s} > 0$;

(v) for any indices $\{t_s \to \infty\}$ if $\limsup_s v_{t_s} < +\infty$ then $\limsup_s |w_{t_s}| < +\infty$;

Then $\lim_{t \to \infty} v_t = 0$, with probability one.